# Poincaré Series of Monomial Rings with Minimal Taylor Resolution

Yohannes Tadesse

October 13, 2011

### Abstract

We give a comparison between the Poincaré series of two monomial rings: $R = A/I$ and $R_q = A/I_q$ where $I_q$ is a monomial ideal generated by the $q$'th power of monomial generators of $I$. We compute the Poincaré series for a new class of monomial ideals with minimal Taylor resolution. We also discuss the structure of a monomial ring with minimal Taylor resolution where the ideal is generated by quadratic monomials.

## 1   Introduction

Let $A = k[x_1, \ldots, x_n]$ be a polynomial ring over a field $k$ and $I$ be an ideal of $A$. The Poincaré series of $R = A/I$ is the power series

$$P^R_k(z) = \sum_{i \geq 0} \dim_k(\operatorname{Tor}_i^R(k,k)) z^i \in \mathbb{Z}[[z]].$$

It is the generating function of the sequence of Betti numbers of a minimal free resolution of $k$ over $R$. The question if this is a rational function was asked by Serre. An affirmative answer was presented by Backelin in [3] when $I$ is a monomial ideal in $A$ and a counter-example was given by Anick in [1,2] when $I = (x_1^2, x_2^2, x_4^2, x_5^2, x_1x_2, x_4x_5, x_1x_3 + x_3x_4 + x_2x_5) \subset A = k[x_1, \ldots, x_5]$.

In this paper we compute the Poincaré series of some monomial rings. Recall that a monomial ring is the quotient ring $R = A/I$ where $I$ is a monomial ideal in the polynomial ring $A$. There exists a unique finite set of monomial generators for a monomial ideal $I$ which we denote by $G(I) = \{m_1, \ldots, m_t\}$. We also denote the graded maximal ideal of $A$ by $\mathfrak{m}_A = (x_1, \ldots, x_n)$, the set containing the first $d$ positive integers by $\mathbb{N}_d$ and by $|S|$ the cardinality of a set $S$.

We will use polarization of monomial rings introduced by Fröberg in [6] to prove that, for any positive integer $q$, if $R = A/I$ is a monomial ring such that $I \subseteq \mathfrak{m}_A^2$ and $R_q = A/I_q$ is a monomial ring such that $G(I_q) = \{m^q \mid m \in G(I)\}$, then $P_k^R(z) = P_k^{R_q}(z)$, (Theorem 2.3). This result does not hold if there exists some $x_i \in G(I)$, that is $I \subsetneqq \mathfrak{m}_A^2$, as the following example shows.

**Example 1.1.** Let $R = A/I$ where $I = (x) \subset A = k[x]$. Then $P_k^R(z) = 1$ but, since $R_2 = A/(x^2)$ is a complete intersection, we have $P_k^{R_2}(z) = \frac{1+z}{1-z^2} = \frac{1}{1-z}$.

We will give the minimal generating set of the homology $H(K^R)$ of the Koszul complex $K^R$ of a certain class of monomial rings $R = A/I$ with minimal Taylor resolution (Theorem 3.3). Namely, we consider monomial ideals $I$ such that there exists a total ordering $m_1 \prec m_2 \prec \cdots \prec m_t$ on $G(I)$ and a positive integer $d \leq t$ such that $\gcd(m_i, \ldots, m_{i+d-1}) \neq 1$ and $\gcd(m_i, m_{i+d+j}) = 1$ for all $i = 1, \ldots, t-d$ and $j > 1$. We use this description to compute the Hilbert series of $H(K^R)$. Then Fröberg in [7] proved that $P_k^R(z)$ is a quotient of the Hilbert series of two associative graded algebras; namely $K^R$ and $H(K^R)$.

It is of some interest to know classes of monomial rings with a minimal Taylor resolution which are either a complete intersection or Golod. The reader can see the definition of a Golod ring in e.g. [7]. It is well known that a monomial ring $R = A/I$ is Golod if $G(I)$ has a common factor $\neq 1$. We call such a ring trivially Golod. Since Poincaré series of such rings are already known, one may determine the structure of $R = A/I$ from $P_k^R(z)$. In section 4 we consider monomial rings $R = A/I$ with minimal Taylor resolution when $I$ is either a stable ideal, an ideal with a linear resolution and an ideal generated by quadratic monomials. We study the conditions on $G(I)$ for such a monomial ideal $I$ so that $R = A/I$ becomes a complete intersection or trivially Golod.

## 2   Ideals generated by Powers of generators

For a monomial $m \in A$, let $\operatorname{Supp}(m) = \{j \mid x_j \text{ divides } m\}$. Given a monomial ideal $I \subset A$, if $j \notin \bigcup_{m \in G(I)} \operatorname{Supp}(m)$, then we have $R = A/I \cong A'/(I \cap A') \otimes_k k[x_j]$ where $A' = k[x_1, \ldots, \hat{x}_j, \ldots, x_n]$. Since $P_k^{k[x_j]}(z) = 1 + z$ and the Poincaré series of a tensor product of two algebras is the product of Poincaré series of the algebras, it follows that $P_k^R(z) = (1 + z)P_k^{A'/(I \cap A')}(z)$. So we may assume, without any loss of generality, that $\bigcup_{m \in G(I)} \operatorname{Supp}(m) = \{1, \ldots, n\}$. Furthermore, let $I \subset A$ be a monomial ideal such that $G(I) = \{x_1, \ldots, x_s, m_1, \ldots, m_t\}$ for some $s \leq n$ where $\deg(m_i) > 1$ for each $i = 1, \ldots, t$. Since then $\cup_{i=1}^t \operatorname{Supp}(m_i) = \{s+1, \ldots, n\}$, it follows that $R \cong k \otimes_k R' \cong R'$ where $R' = k[x_{s+1}, \ldots, x_n]/(m_1, \ldots, m_t)$. Hence $P_k^R(z) = P_k^{R'}(z)$. So we shall consider monomial rings $R = A/I$ such that $\deg(m) > 1$ for each $m \in G(I)$, that is $I \subseteq \mathfrak{m}_A^2$.

The following Proposition is due to Fröberg in [6], pp. 30-32. We include it for reference since it plays an important role in proving Theorem 2.3.

**Proposition 2.1.** *(Polarization) Let $R = \boldsymbol{k}[x_1, \ldots, x_n]/I$ be a monomial ring. There exists a squarefree monomial ring $S = \boldsymbol{k}[y_1, \ldots, y_N]/I'$ such that $R = S/(f_1, \ldots, f_{N-n})$, where $f_1, \ldots, f_{N-n}$ is a regular sequence of homogeneous elements of degree one. Moreover,*

   *1. $R$ is complete intersection if and only if $S$ is complete intersection.*

   *2. $R$ is Golod if and only if $S$ is Golod.*

   *3. $P_{\boldsymbol{k}}^R(z) = \frac{P_{\boldsymbol{k}}^S(z)}{(1+z)^{N-n}}$.*

Recall that a squarefree monomial ring is the quotient of a polynomial ring by a squarefree monomial ideal. We shall explain how to get the squarefree monomial ring $S$ from $R$. Let

$$(2.1) \qquad I(x_i) = \max\{\alpha \in \mathbb{Z}_{\geq 0} \mid x_i^\alpha \text{ divides some } m \in G(I)\}$$

If each $I(x_i) = 1$, then $R$ is squarefree. If there exists some $i$ with $I(x_i) > 1$, we introduce new variables and replace each monomial $m \in G(I)$ by a squarefree monomial of degree equal to $\deg(m)$ in $N = \sum_{i=1}^n I(x_i)$ new variables, say $y_1, \ldots, y_N$. This set of squarefree monomials generate a monomial ideal $I'$ in the polynomial ring $A' = k[y_1, \ldots, y_N]$ and we take $S = A'/I'$.

**Example 2.2.** Consider the monomial ring $R = k[x_1, x_2, x_3]/(x_1^3, x_2^2 x_3, x_1 x_2 x_3)$. Then $N = 3+2+1 = 6$. Replacing $x_1^3$ by $y_1 y_2 y_3$, $x_2^2 x_3$ by $y_4 y_5 y_6$ and $x_1 x_2 x_3$ by $y_1 y_4 y_6$ we obtain a squarefree monomial ideal $I' = (y_1 y_2 y_3, y_4 y_5 y_6, y_1 y_4 y_6) \subset A' = k[y_1, \ldots, y_6]$. It is not difficult to see that $R \cong S/(y_1-y_2, y_1-y_3, y_4-y_5)$ where $S = A'/I'$.

**Theorem 2.3.** *Let $R = A/I$ be a monomial ring such that $I \subseteq \mathfrak{m}_A^2$, and $R_q = A/I_q$ be a monomial ring such that $G(I_q) = \{m^q \mid m \in G(I)\}$ for some integer $q > 1$. Then $P_k^R(z) = P_k^{R_q}(z)$.*

*Proof.* We consider two cases depending on the values of $I(x_i)$ defined in (2.1).

(a): Let $I(x_i) = 1$ for all $i \in \mathbb{N}_n$. Then $I$ is a squarefree monomial ideal. It suffices to prove that $\operatorname{Tor}_i^{R_q}(k, k) \cong \operatorname{Tor}_i^R(k, k)$ for all $i \geq 0$. Put $A_q = k[x_1^q, \ldots, x_n^q]$ and $R' = A_q/A_q I_q$. There exists a natural isomorphism of $k$-algebras $A \to A_q$ and, since $I$ is a squarefree monomial ideal, this induces an isomorphism of $k$-algebras $R \to R'$. On the other hand, since $A_q I_q \subset I_q$ we have a ring map $R' \to R_q$ which defines an $R'$-module structure on $R_q$.

In fact one has $R_q = \oplus_{|\alpha| < q} R'(X^\alpha \mod I_q)$, where $\alpha \in \mathbb{Z}_{\geq 0}^n$, so $R_q$ is free over $R'$. Considering $k$ both as an $R'$-module and an $R_q$-module, and using the ring map $R' \to R_q$, it follows, by [9, Prop. 3.2.9], that

$$\mathrm{Tor}_i^{R'}(k, k) \cong \mathrm{Tor}_i^{R_q}(k \otimes_{R'} R_q, k)$$

for all $i \geq 0$. Using the projection map $R_q \to R' \cong R$, we obtain

$$k \hookrightarrow k \otimes_R R_q \to k \otimes_R R \cong k.$$

So we have $\mathrm{Tor}_i^R(k \otimes_{R_q} R, k) \cong \mathrm{Tor}_i^R(k, k)$.

(b): Let $I(x_i) > 1$ for some $i$. We use Prop. 2.1 and (a) above to prove $P_k^R(z) = P_k^{R_q}(z)$. Put $N = \sum_{i=1}^n I(x_i)$. By Prop. 2.1 there exists a square free monomial ideal $J \subset B = k[y_1, \ldots, y_N]$ with $G(J) = \{M_1, \ldots, M_t\}$ and a regular sequence $f_1, \ldots, f_{N-n} \in S = B/J$ of degree one such that $R \cong S/(f_1, \ldots, f_{N-n})$. It also follows by Prop. 2.1 that

$$(2.2) \qquad P_k^R(z) = \frac{P_k^S(z)}{(1+z)^{N-n}}$$

Now for an integer $q > 1$, let $J_q \subset B$ be a monomial ideal with $G(J_q) = \{M_1^q, \cdots, M_t^q\}$ and put $S_q = B/J_q$. By (a) above we have

$$(2.3) \qquad P_k^S(z) = P_k^{S_q}(z).$$

Since each $M_j$ is a squarefree monomial, $j \in \mathbb{N}_t$, we have $J_q(y_i) = q > 1$ for each $i \in \mathbb{N}_N$. Again by Prop. 2.1 there exists a square free monomial ideal $J' \subset C = k[z_1, \ldots, z_{Nq}]$ and a regular sequence $g_1, \ldots, g_{Nq-N} \in S' = C/J'$ such that $S_q \cong S'/(g_1, \ldots, g_{Nq-N})$ and, moreover,

$$(2.4) \qquad P_k^{S_q}(z) = \frac{P_k^{S'}(z)}{(1+z)^{Nq-N}}$$

Combining (2.2-2.4), we obtain $P_k^R(z) = \frac{P_k^{S'}(z)}{(1+z)^{Nq-n}}$.

On the other hand, since $G(I_q) = \{m_1^q, \ldots, m_t^q\}$, one has $I_q(x_i) = q \cdot I(x_i) > 1$ for each $i$. So there exists a square free monomial ideal $I'_q \in C = k[z_1, \ldots, z_{Nq}]$ and a regular sequence $h_1, \ldots, h_{Nq-n} \in S' = C/I'_q$ such that $R_q \cong S'/(h_1, \ldots, h_{Nq-n})$. So $P_k^{R_q}(z) = \frac{P_k^{S'}(z)}{(1+z)^{Nq-n}}$. We thus have $P_k^{R_q}(z) = \frac{P_k^{S'}(z)}{(1+z)^{Nq-n}} = P_k^R(z)$. $\qquad \square$

## 3 Rings with a Minimal Taylor Resolution

Let $R = A/I$ be a monomial ring with $G(I) = \{m_1, \ldots, m_t\}$ and $T$ be the exterior algebra of a rank $t$ free $A$-module with a standard basis $e_{i_1, \ldots, i_l}$ for

$1 \leq i_1 < \cdots < i_l \leq t$. Consider $T$ as a finite free resolution of $R$ with a differential

$$d(e_{i_1,\ldots,i_l}) = \sum_{j=1}^{l} (-1)^{j-1} \frac{m_{i_1,\ldots,i_l}}{m_{i_1,\ldots,\hat{i}_j,\ldots,i_l}} e_{i_1,\ldots,\hat{i}_j,\ldots,i_l}$$

where $m_{i_1,\ldots,i_l} = \text{lcm}(m_{i_1},\ldots,m_{i_l})$. This resolution is called the Taylor resolution of $R$, see also [5]. It is far from being minimal. But we get a minimal Taylor resolution whenever $m_{i_1,\ldots,i_l} \neq m_{i_1,\ldots,\hat{i}_j,\ldots,i_l}$ for all $i_1,\ldots,i_l \in \mathbb{N}_t$, or equivalently, whenever each $m_i$ contains a variable with a maximal power. The reader may refer [7] for other equivalent conditions. It is evident that if a monomial ring $R = A/I$ has a minimal Taylor resolution, then so will the monomial ring $R_q = A/I_q$ where $G(I_q) = \{m^q \mid m \in I\}$ for an integer $q > 0$.

Let $R = A/I$ be a monomial ring. A differential graded, associative and commutative algebra structure for the Taylor resolution of $R$ was given by Gemeda in [4]. Using this algebra structure, Fröberg in [7, Theorem 3] proved that the Poincaré series of $R = A/I$ having a minimal Taylor resolution is the quotient of the Hilbert series of the Koszul complex $K^R$ of $R$ and the Hilbert series of its Homology $H(K^R)$. More precisely,

$$P_k^R(z) = \frac{Hilb(K^R)(z)}{Hilb(H(K^R))(-z,z)} = \frac{(1+z)^n}{Hilb(H(K^R))(-z,z)}$$

where $n$ is the embedding dimension of $R$ and we consider $H(K^R)$ as a bi-graded algebra by a polynomial degree and a total degree. The basis for $H(K^R)$ can be described in terms of representatives $T_1,\ldots,T_n$ in $K^R$ by

$$f_{i_1,\ldots,i_l} = \frac{\text{lcm}(m_{i_1},\ldots,m_{i_l})}{x_{i_1}\cdots x_{i_l}} T_{i_1} \cdots T_{i_l}$$

for $1 \leq i_1 < \ldots < i_l \leq n$.

**Example 3.1.** We will describe how to compute $P_k^R(z)$ for a monomial ring $R = A/I$ where $I = (x^2y, y^2z, z^2) \subset A = k[x,y,z]$. Let $T_1, T_2, T_3$ be the standard generators for $K^R$. Then

$$f_1 = \frac{x^2y}{x}T_1 = xyT_1 \qquad\qquad f_2 = \frac{y^2z}{y}T_2 = yzT_2$$

$$f_3 = \frac{z^2}{z}T_3 = zT_3 \qquad\qquad f_{12} = \frac{lcm(x^2y,y^2z)}{xy}T_1T_2 = xyzT_1T_2$$

$$f_{23} = \frac{lcm(y^2z,z^2)}{yz}T_2T_3 = yzT_2T_3 \qquad f_{13} = \frac{lcm(x^2y,z^2)}{xz}T_1T_3 = xyzT_1T_3$$

$$f_{123} = \frac{lcm(x^2y,y^2z,z^2)}{xyz} = xyzT_1T_2T_3$$

We then have $\bar{f}_1 \bar{f}_3 = \bar{f}_{13}$ and, otherwise, $\bar{f}_I \bar{f}_J = 0$ for all $I, J \subset \{1, 2, 3\}$. So we obtain $H(K^R) = k(X_1, X_2, X_3, X_{12}, X_{23}, X_{123}/(X_I X_J \mid$ for all $I, J$ except $I = \{1\}$, and $J = \{3\}$)). The Hilbert series, then, becomes

$$Hilb(H(K^R))(X, Y) = 1 + 3XY + 2XY^2 + X^2Y^2 + XY^3$$

and the Poincaré series is

$$P_k^R(z) = \frac{Hilb(K^R)(z)}{Hilb(H(K^R))(-z, z)} = \frac{(1+z)^3}{1 - 3z^2 - 2z^3}.$$

A *strictly ordered partition* of a finite totally ordered set $(S, \prec)$ is a sequence $(S_1, \ldots, S_l)$ of non-empty subsets of $S$ such that they form an ordered partition and $\max(S_i) \prec \min(S_{i+1})$ for all $i = 1, \ldots, l-1$. In this case we call $l$ the *length* and $(|S_1|, \ldots, |S_l|)$ the *weight* of the partition. The following is evident.

**Proposition 3.2.** *Fix positive integers $d \leq t$. For any non-empty subset $S$ of $\mathbb{N}_t$ there exists a strictly ordered partition $(S_1, \ldots, S_l)$ such that:*

1. *Any two consecutive numbers in $S_i$ differ at most by $d-1$.*

2. *$\min(S_{j+1}) - \max(S_j) \geq d$ for each $j = 1, \ldots, l$.*

**Theorem 3.3.** *Fix a positive integer $d$. Let $R = A/I$ be a monomial ring with a minimal Taylor resolution. Assume that $|G(I)| = t$ and there exists a total ordering on $G(I)$ such that $\gcd(m_i, m_{i+1}, \ldots, m_{i+d}) \neq 1$ and $\gcd(m_i, m_{i+d+j}) = 1$ for any $j \geq 0$ and $i = 1, \ldots, t-d$. Consider the collection*

$$\mathcal{B} = \{S \subseteq \mathbb{N}_t \mid \text{any two consecutive numbers in } S \text{ differ at most by } d-1\}.$$

*We have the following:*

1. *The minimal generating set of the homology $H(K^R)$ of $K^R$ is $\{X_S \mid S \in \mathcal{B}\}$, that is $H(K^R) = k(X_S)_{S \in \mathcal{B}}/J$ where $J$ is the ideal*

   $J = (X_{S_1} \cdots X_{S_l} \mid$ each $S_i \in \mathcal{B}, \quad S_i \cup S_j \notin \mathcal{B}$ for all $i, j \in \mathbb{N}_l$ and
   $S_i \cap S_j \neq \emptyset$ for some $i, j \in \mathbb{N}_l$).

2. *For any non-empty set $S \subseteq \mathbb{N}_t$ there exists a partition $S_1, \ldots, S_l$ of $S$ such that each $S_j \in \mathcal{B}$ and $S_i \cup S_j \notin \mathcal{B}$. Put $m = |S|$. Then the Hilbert series of $H(K^R)$ is*

   $$(3.1) \qquad Hilb(H(R^K))(X, Y) = \sum_{(n_1, \ldots, n_l)} f(n_1, \ldots, n_l) X^l Y^m.$$

   *where $f(n_1, \ldots, n_l)$ is the number of subsets of $N_t$ having a strictly ordered partition defined in Prop. 3.2 for a given length $l \leq t$ and a weight $(n_1, \ldots, n_l) \in \mathbb{Z}_{>0}^l$.*

*Proof.* (1): Note that $H(K^R) = k(X_S)_{S \subseteq \mathbb{N}_t}/J$. Now let $(S_1, \ldots, S_l)$ be a strictly ordered partition given in Prop. 3.2 of a set $S \subseteq \mathbb{N}_t$. Then each $S_i \in \mathcal{B}$ and by assumption any two monomials in $G(I)$ indexed by elements of different subpartitions are relatively prime. We then have $\bar{f}_S = \bar{f}_{S_1} \cdots \bar{f}_{S_l}$ and so $X_S = X_{S_1} \cdots X_{S_l}$. It follows that set $\{X_{S_i} \mid S_i \in \mathcal{B}\}$ generates $H(K^R)$. Since each $S_i \in \mathcal{B}$ has a strictly ordered partition of length 1, $\{X_{S_i} \mid S_i \in \mathcal{B}\}$ becomes a minimal generating set. (2): Follows from Prop 3.2. $\qquad\square$

Now we give an example.

**Example 3.4.** Let $R = A/I$ where $I = (x^2yz, y^2zw, z^2wu, w^2u, u^2) \subset A = k[x, y, z, w, u]$. We want to compute $Hilb(H(K^R))(X, Y)$. Consider the ordering $m_1 = x^2yz \prec m_2 = y^2zw \prec m_3 = z^2wu \prec m_4 = w^2u \prec m_5 = u^2$ where $d = 3$. Then $\mathcal{B}$ consists of all non-empty subsets of $\mathbb{N}_5$ except $\{1, 4\}, \{1, 5\}, \{2, 5\}, \{1, 2, 5\}$ and $\{1, 4, 5\}$. That is, $\bar{f}_{14} = \bar{f}_1 \bar{f}_4$, $\bar{f}_{15} = \bar{f}_1 \bar{f}_5$, $\bar{f}_{25} = \bar{f}_2 \bar{f}_5$, $\bar{f}_{125} = \bar{f}_{12} \bar{f}_5$ and $\bar{f}_{145} = \bar{f}_1 \bar{f}_{45}$. Therefore,

$$Hilb(H(K^R))(X, Y) = 1 + 5XY + 8XY^2 + 2X^2Y^2 + 8XY^3 + 2X^2Y^3 + 5XY^4 + XY^5.$$

We obtain a formula for Hilbert series of $H(K^R)$ if $d \leq 2$ in Theorem 3.3.

**Proposition 3.5.** *Keep all the assumptions of Theorem 3.3 for a monomial ring $R = A/I$.*

1. *If $d = 1$, $R$ is a complete intersection. If $d = t$, $R$ is trivially Golod.*

2. *If $d = 2$, the Hilbert series of $H(K^R)$ is*

$$(3.2) \qquad Hilb(H(R^K))(X, Y) = \sum_{(n_1, \ldots, n_l)} \binom{t - m + 1}{l} X^l Y^m.$$

*Proof.* (1) is clear, so we prove only (2). If $d = 2$, the strictly ordered partition of a non-empty subset $S$ of $\mathbb{N}_t$ is given by a partition $(S_1, \ldots, S_l)$ such that each $S_i$ contains consecutive integers and $\min(S_i) - \max(S_{i-1}) \geq 2$. Now let $S_i = \{a_i, a_i + 1, \ldots, a_i + n_i - 1\}$ where $a_i = \min(S_i)$ for each $i$. We then obtain the following inequalities:

$$1 \leq a_1$$
$$a_1 + n_1 - 1 + 2 \leq a_2 \Rightarrow a_1 + n_1 < a_2$$
$$a_2 + n_2 - 1 + 2 \leq a_3 \Rightarrow a_1 + n_1 + n_2 < a_3$$
$$a_3 + n_3 - 1 + 2 \leq a_4 \Rightarrow a_1 + n_1 + n_2 + n_3 < a_4$$
$$\vdots$$
$$a_{l-1} + n_{l-1} - 1 + 2 \leq a_l \Rightarrow a_1 + \sum_{i=1}^{l-1} n_i < a_l$$
$$a_l + n_l - 1 \leq t.$$

This is equivalent to the inequality system $1 \leq a_1 < a_2 - n_1 < a_3 - (n_1 + n_2) < \cdots < a_l - (\sum_{i=1}^{l-1} n_i) \leq t - m + 1$. The number of solutions we get for this inequality is $\binom{t-m+1}{l}$. $\qquad \square$

**Remark 3.6.** It is known that for a monomial ring $R = A/I$ with minimal Taylor resolution, the dimension of the $m$'th homology of $K^R$ is $\binom{|G(I)|}{m}$, see [7]. This value also equals to the sum of the coefficients in $Hilb(H(R^R))(X,Y)$ with terms containing $Y^m$. It then follows from Prop. 3.5 that

$$\sum_{\substack{(n_1,\ldots,n_l) \\ \sum_i n_i = m}} \binom{|G(I)| - m + 1}{l} = \binom{|G(I)|}{m}$$

where $(n_1, \ldots, n_l)$ is the weight of a strictly ordered partition defined in Prop. 3.2 for a set $S \subset \mathbb{N}_t$ with $|S| = m$ and $d = 2$.

# 4    Complete Intersection and Trivially Golod Rings

From the algorithm to compute $P_k^R(z)$ given in [7] it follows, for a monomial ring $R = A/I$ with minimal Taylor resolution, that $P_k^R(z) = (1-z)^n/(1-z^2)^t$ if and only if $\gcd(m_i.m_j) = 1$ for all $i \neq j$ where $G(I) = \{m_1, \ldots, m_t\}$, i.e. $R$ is a complete intersection. Furthermore $P_k^R(z) = (1+z)^n/(1 - \sum_{i=1}^{t} \binom{t}{i} z^{i+1})$ if $G(I)$ has a common factor $\neq 1$, i.e. $R$ is trivially Golod. This gives Prop. 4.1.

**Proposition 4.1.** *Let $R = A/I$ and $R_q = A/I_q$ be two monomial rings such that $G(I_q) = \{m^q \mid m \in G(I)\}$ for some integer $q > 0$. Then*

1. *$R_q$ is trivially Golod if and only if $R$ is trivially Golod.*

2. *$R_q$ is complete intersection if and only if $R$ is complete intersection.*

For a monomial $m$, put $i_0 = \max(\text{Supp}(m))$. Recall that a monomial ideal $I$ is said to be stable if for each monomial $m \in I$ and all $i < i_0$, we have $x_i m / x_{i_0} \in I$. In [8] Okudaira and Takayama proved that such an ideal $I$ has a minimal Taylor resolution if and only if each $m_i \in G(I)$ has the form $m_i = x_i (\prod_{j=1}^{i} x_j^{n_j})$ for $i = 1, \ldots, t$ and for some integers $n_1, \ldots n_r \geq 0$. It follows that $R$ is trivially Golod if $n_1 > 0$; and $R$ is a complete intersection if each $n_i = 0$.

A monomial ideal $I$ with a linear resolution possesses a minimal Taylor resolution if and only if each $m_i \in G(I)$ is of the form $m_i = u x_{j_i}$ for some $j_i \in \mathbb{N}_n$ and a monomial $u \in A$, see [8]. It follows then that $R$ is a complete intersection if $u = 1$, it is trivially Golod if $u \neq 1$.

**Theorem 4.2.** *Let $R = A/I$ be a monomial ring with a minimal Taylor resolution and each $m \in G(I)$ is a quadratic monomial. Then $R$ is a $k$-tensor product of a complete intersection and trivially Golod rings.*

*Proof.* Let $P_1, \ldots, P_r$ be a partition of $G(I)$ such that any two elements of one subpartition have a common factor and elements between any pair of different subpartitions are relatively prime. Put $A_i = k[x]_{x \in \text{Supp}(P_i)}$. Since $\text{Supp}(G(I))$ is partitioned by the collection $\{\text{Supp}(P_i)\}_i$, we have $A = \otimes_{i=1}^{r} A_i$. If each $P_i$ is a singleton, by construction, $R$ is a complete intersection. Since the monomials in each partitions are quadratic, there exists a $j \in \mathbb{N}_n$ such that $\gcd(m)_{m \in P_i} = x_j$. So we obtain a monomial ideal $I_i = (m)_{m \in P_i} \subset A_i$ such that $I = \sum_i I_i$, $R = \otimes A_i / I_i$ and each $R_i = A_i / I_i$ has a minimal Taylor resolution. Moreover, if $I_i$ is principal, then $R_i$ is a complete intersection and otherwise $R_i$ is trivially Golod. $\qquad\square$

**Example 4.3.** Consider the monomial ideal $I = (x_1^2, x_1 x_2, x_1 x_4, x_3^2, x_5^2) \subset A = k[x_1, \ldots, x_5]$. It is easy to see that $R = A/I$ has a minimal Taylor resolution. We have a partition of $G(I)$ in three parts $P_1 = \{x_1^2, x_1 x_2, x_1 x_4\}$, $P_2 = \{x_3^2\}$, $P_3 = \{x_5^2\}$ and monomial ideals $I_1 = (x_1^2, x_1 x_2, x_1 x_4) \subset A_1 = k[x_1, x_2, x_4]$, $I_2 = (x_3^2) \subset A_2 = k[x_3]$ and $I_3 = (x_5^2) \subset A_3 = k[x_5]$. We thus have a trivially Golod ring $R_1 = A_1 / I_1$, complete intersections $R_2 = A_2 / I_2$ and $R_3 = A_3 / I_3$. So $R = R_1 \otimes_k R'$ where $R' = R_2 \otimes_k R_3$ is a complete intersection.

## Acknowledgment

The author would like to thank Prof. Ralf Fröberg for presenting the problem in Pragmatic 2011, for the discussion we have had and for his valuable suggestion. The author also would like to thank Jörgen Backelin for his suggestion in proving Prop. 3.5.

Department of Mathematics,
Stockholm University,
SE 106-91,
Stockholm, SWEDEN,
*E-mail address: tadesse@math.su.se*